\documentclass[11pt,leqno,fleqn]{article}
\usepackage{epsf,amsmath,amssymb,amscd,hyperref,srcltx}
\usepackage{amsfonts}
\usepackage{counters}
\usepackage{psfrag}
\newcommand{\dps}{\displaystyle}

\newcommand{\T}{^{\sf T}}
\newcounter{figuren}

\newcommand{\dy}[2]{%
\refstepcounter{equation}%
\LABEL{#1}%
\begin{list}{}{
\topsep 5mm
\leftmargin 18mm
\rightmargin 0cm
\itemsep 0mm
\listparindent 0mm
\parsep 0mm
\itemsep 0mm
\labelsep 0mm
\labelwidth 18mm
}%
\item[\rm (\theequation)\hfill]
#2
\end{list}%
}
\newcommand{\dyz}[1]{%
\refstepcounter{equation}%
\begin{list}{}{
\topsep 5mm
\leftmargin 18mm
\rightmargin 0cm
\itemsep 0mm
\listparindent 0mm
\parsep 0mm
\itemsep 0mm
\labelsep 0mm
\labelwidth 18mm
}%
\item[\rm (\theequation)\hfill]
#1
\end{list}%
}
\newcommand{\dyyz}[1]{\dyz{\raggedright$\dps#1$}}
\newcommand{\dyy}[2]{\dy{#1}{\raggedright$\dps#2$}}
\newcommand{\de}[2]{\dy{#1}{\raggedright$\displaystyle #2 $}}
\newcommand{\dez}[1]{\dyz{\raggedright$\displaystyle #1 $}}
\newcommand{\leeg}[1]{}

\newcounter{stelling}
\newcommand{\thm}[2]{\setcounter{gevolg}{0}\setcounter{claim}{0}\refstepcounter{stelling}\vspace{4mm}\noindent{\bf Theorem \thestelling.}\label{#1}{\it #2}}

\newcommand{\thmz}[1]{\setcounter{gevolg}{0}\setcounter{claim}{0}\refstepcounter{stelling}\vspace{4mm}\noindent{\bf Theorem \thestelling.}{\it #1}}

\newcommand{\thmnmz}[2]{\setcounter{gevolg}{0}\setcounter{claim}{0}\refstepcounter{stelling}\vspace{4mm}\noindent{\bf Theorem \thestelling\ {\rm (#1)}.}{\it #2}}

\newalphacounterwithin{gevolg}{stelling}

\newcounter{hulpstelling}
\newcommand{\lemma}[2]{\refstepcounter{hulpstelling}\setcounter{claim}{0}\vspace{4mm}\noindent{\bf Lemma \thehulpstelling.}\label{#1}{\it #2}}

\newcounter{bewering}

\newcounter{claim}

\newalphacounterwithin{subclaim}{claim}

\newcounter{opmerking}

\newcounter{hoofdstuk}

\newcounter{sectie}
\newcounter{subsectie}
\newcounterwithin{ex}{sectie}

\newcommand{\sectz}[1]{\refstepcounter{sectie}\setcounter{subsectie}{0}\setcounter{ex}{0}
\section*{\boldmath \thesectie. #1}%
}

\newcounter{lit}

\newcommand{\pf}{\vspace{3mm}\noindent{\bf Proof.}\ }

\newcommand{\bx}{\hspace*{\fill} \hbox{\hskip 1pt \vrule width 4pt height 8pt depth 1.5pt \hskip 1pt}

\addvspace{4mm}}

\newcommand{\rf}[1]{{\rm (\ref{#1})}}

\newcommand{\kint}[2]{\mbox{$\int$}}

\newcommand{\NIET}[1]{}

\newcommand{\LABEL}[1]{\label{#1}}

\newcommand{\oN}{{\mathbb{N}}}

\newcommand{\oR}{{\mathbb{R}}}

\begin{document}

\begin{center}
{\LARGE\bf Weak and strong regularity, compactness, and approximation of polynomials

}
\vspace{4mm}

{\large
\hspace{10mm}
Alexander Schrijver\footnote{ CWI and University of Amsterdam.
Mailing address: CWI, Science Park 123, 1098 XG Amsterdam,
The Netherlands.
Email: lex@cwi.nl.}}

\end{center}

\noindent
{\small{\bf Abstract.}
Let $X$ be an inner product space,
let $G$ be a group of orthogonal transformations of $X$,
and let $R$ be a bounded $G$-stable subset of $X$.
We define very weak and very strong regularity for such pairs $(R,G)$
(in the sense of Szemer\'edi's regularity lemma), and prove that
these two properties are equivalent.

Moreover, these properties are
equivalent to the compactness of the space $(B(H),d_R)/G$.
Here $H$ is the completion of $X$ (a Hilbert space),
$B(H)$ is the unit ball in $H$,
$d_R$ is the metric on $H$ given by
$d_R(x,y):=\sup_{r\in R}|\langle r,x-y\rangle|$,
and $(B(H),d_R)/G$ is the orbit space of $(B(H),d_R)$
(the quotient topological space with the $G$-orbits as quotient classes).

As applications we give Szemer\'edi's regularity lemma,
a related regularity lemma for partitions into intervals,
and a low rank approximation theorem for homogeneous polynomials.

}

\sectz{Equivalence of very weak regularity, very strong regularity, and
compactness}

This paper is inspired by Szemer\'edi's regularity lemma
([7]) and subsequent work on graph limits
by Lov\'asz and Szegedy ([3,\linebreak[0]4])
(cf.\ also [5]).

Let $X$ be an inner product space and let $R$ be a bounded subset of $X$
spanning $X$.
(So each element of $X$ is a linear combination of finitely many elements
of $R$.)
Let $G$ be a group of orthogonal transformations $\pi$ of $X$
with $\pi(R)=R$.
Let $B(X)$ denote the unit ball in $X$.
For any $k$, let $R_k:=\{\pm r_1\pm\cdots\pm r_k\mid
r_1,\ldots,r_k\in R\}$.
Let $H$ be the completion of $X$, which is a Hilbert space.
Then $G$ naturally acts on $H$.
For $x,y\in H$, define
\dez{
d_R(x,y):=\sup_{r\in R}|\langle r,x-y\rangle|.
}
The space $(B(H),d_R)/G$ is the orbit space of $(B(H),d_R)$, i.e., the
quotient topological space of $(B(H),d_R)$ taking the $G$-orbits as classes.

\thm{3ok12b}{
The following are equivalent:\\
{\rm (i)} $(R,G)$ is \textit{\textbf{very weakly regular}}:
for each $k$ there exists a finite set $Z\subseteq X$
such that for each $x\in R_k$ there exist $z\in Z$ and $\pi\in G$
satisfying $\langle r,x-z^{\pi}\rangle^2\leq 1$ for each $r\in R$;\\
{\rm (ii)} $(R,G)$ is \textit{\textbf{weakly regular}}:
for each $\varepsilon>0$ there exists a finite set $Z\subseteq B(X)$
such that for each $x\in B(X)$ there exist $z\in Z$ and $\pi\in G$
satisfying $|\langle r,x-z^{\pi}\rangle|<\varepsilon$
for each $r\in R$;
{\rm (iii)} $(R,G)$ is \textit{\textbf{very strongly regular}}:
for each $\varepsilon>0$ and $f:X\to\{1,2,\ldots\}$
there exists a finite set $Z\subseteq B(X)$
such that for each $x\in B(X)$
there exist $z\in Z$ and $\pi\in G$ satisfying\footnote{$\|.\|_p$ is the $L^p$-norm, here for
the finite-dimensional space $\oR^{f(z)}$.}
\dyy{25se12a}{
\sum_{i=1}^{f(z)}
|\langle r_i,x-z^{\pi}\rangle|^t
\leq
\varepsilon
(1+\|(\|r_1\|^t,\ldots,\|r_{f(z)}\|^t)\|_p),
}
for all $t\in[\varepsilon,2]$, where $p:=2/(2-t)$,
and for all orthogonal $r_1,\ldots,r_{f(z)}\in R$;\\
{\rm (iv)} the space $(B(H),d_R)/G$ is compact.
}

\pf
(iii)$\Rightarrow$(ii) follows by taking $f(x)=1$ for each $x\in X$
and $t=1$.
(ii)$\Rightarrow$(i) follows by observing that $\frac{1}{t}R_k\subseteq B(X)$
for some $t$, and taking $\varepsilon:=1/t$.
So it suffices to prove
(i)$\Rightarrow$(ii),
(ii)$\Rightarrow$(iv),
and
(iv)$\Rightarrow$(iii).

For all $x,y\in B(H)$ define
\dez{
\delta_R(x,y):=\inf_{\pi\in G}d_R(x,y^{\pi}).
}
Then $\delta_R$ is a pseudometric, and the space $(B(H),\delta_R)$ is topologically
homeomorphic to the orbit space $(B(H),d_R)/G$.

Observe that (i) implies that the space $(R_k,\delta_R)$ is totally bounded\footnote{A pseudometric
space is {\em totally bounded} if for each $\varepsilon>0$ it can be covered by finitely
many balls of radius $\varepsilon$ (cf.\ [1]).}.
Indeed, choose $\varepsilon>0$.
Let $t:=\lceil\varepsilon^{-1}\rceil$.
Then $R_{kt}$ can be covered by finitely many $\delta_R$-balls of radius 1.
As $R_k\subseteq\frac{1}{t}R_{kt}$, $R_k$ can be covered by finitely many $\delta_R$-balls
of radius $1/t\leq\varepsilon$.

So we can assume, by scaling, that $\|r\|\leq 1$ for each $r\in R$.

\medskip
\noindent
(i)$\Rightarrow$(ii):
We saw above that (i) implies that $(R_k,\delta_R)$ is totally bounded for each $k$.
Now define, for each $k$,
\dez{
S_k:=\{\lambda_1 r_1+\cdots+\lambda_k r_k\mid r_1,\ldots,r_k\in R, \lambda_1,\ldots,\lambda_k\in[-1,+1]\}.
}
Then also $(S_k,\delta_R)$ is totally bounded.
Indeed, choose $\varepsilon>0$, and define $t:=k\lceil\varepsilon^{-1}\rceil$.
Then each $x\in S_k$ has Hilbert distance less than $\varepsilon$ to $\frac{1}{t}R_{kt}$.
By the above, $(R_{kt},\delta_R)$ is totally bounded, hence so is $(\frac{1}{t}R_{kt},\delta_R)$.
So $(S_k,\delta_R)$ is totally bounded.

Next we show that for each $k$:
\de{3ok12a}{
B(X)\subseteq B_{d_R}(S_k,1/\sqrt{k}).
}
To see this, choose $a\in B(X)$.
Let $a_0:=a$.
If $a_i$ has been found, and $d_R(a_i,0)> 1/\sqrt{k}$,
choose $r$ with $\langle r,a_i\rangle > 1/\sqrt{k}$.
Let $a_{i+1}:=a_i-\langle r,a_i\rangle r$.
Then by induction on $i$, as $\|r\|\leq 1$,
\dyyz{
\|a_{i+1}\|^2=
\|a_i\|^2-2\langle r,a_i\rangle^2+
\langle r,a_i\rangle^2\|r\|^2
\leq
\|a_i\|^2-\langle r,a_i\rangle^2
\leq
\|a_i\|^2-1/k
\leq
1-i/k-1/k
=
1-(i+1)/k.
}
So the process terminates for some $i\leq k$, and we have \rf{3ok12a},
since $a-a_i\in S_k$ and hence
$d_R(a,S_k)\leq d_R(a,a-a_i)=d_R(a_i,0)\leq 1/\sqrt{k}$.

As each $(S_k,\delta_R)$ is totally bounded, \rf{3ok12a} implies that $(B(X),\delta_R)$ is totally bounded.

\medskip
\noindent
(ii)$\Rightarrow$(iv):
By (ii), the space $(B(H),\delta_R)$ is totally bounded.
So it suffices to show that $(B(H),\delta_R)$ is complete.
Let $x_1,x_2,\ldots$ be a Cauchy sequence in $(B(H),\delta_R)$.
We show that it is convergent.
We can assume that $\delta_R(x_n,x_{n+1})<2^{-n}$ for each $n$.
Let $\pi_1$ be the identity in $G$.
For each $n\geq 1$, we can choose $\pi_{n+1}\in G$ such that $d_R(x_n^{\pi_n},x_{n+1}^{\pi_{n+1}})<2^{-n}$.
Replacing $x_n$ by $x_n^{\pi_n}$, we can assume that $x_1,x_2,\ldots$ is a Cauchy sequence in $(B(H),d_R)$.
As $B(H)$ is weakly compact, $x_1,x_2,\ldots$ has a subsequence that converges to some $a\in B(H)$ in the weak topology on $B(H)$.
Then $\lim_{n\to\infty}d_R(x_n,a)=0$.
Indeed, $d_R(x_n,a)\leq 2^{-n+2}$ for each $n$.
Otherwise, $|\langle r,x_n-a\rangle|>2^{-n+2}$ for some $r\in R$.
As $a$ is weak limit of some subsequence of $x_1,x_2,\ldots$, there is an $m\geq n$ with
$|\langle r,x_m-a\rangle|<2^{-n+1}$.
As $|\langle r,x_n-x_m\rangle|\leq d_R(x_n,x_m)<2^{-n+1}$, this gives a contradiction.

\medskip
\noindent
(iv)$\Rightarrow$(iii):
Choose $\varepsilon>0$ and $f:X\to\{1,2,\ldots\}$.
For any $k$,
consider the function $\phi_k:X\to\oR$ defined by
\dyy{5ok12a}{
\phi_k(x):=\sup_{t\in[\varepsilon,2]}
\sup_{\text{orthogonal}\atop r_1,\ldots,r_k\in R}
\frac{
\sum_{i=1}^k
|\langle r_i,x\rangle|^t
}
{
(1+\|(\|r_1\|^t,\ldots,\|r_k\|^t)\|_p)
}
}
for $x\in X$, where $p=(1-t/2)^{-1}$.
Then $\phi_k$ is continuous with respect to the $d_R$-topology on $B(H)$.
To see this, let $x,y\in B(H)$ with $d_R(x,y)\leq 1$.
Then
$|\langle r,x\rangle|^t-|\langle r,y\rangle|^t\leq
2|\langle r,x-y\rangle|^{\varepsilon}\leq
2d_R(x,y)^{\varepsilon}$ for
each $r\in R$ and $t\in[\varepsilon,2]$.\footnote{This follows from the fact that if
$0\leq b\leq a\leq 1$, then for each $t\in[1,2]$:
$a^t-b^t\leq a^t-b^t+(a^{2-t}-b^{2-t})(ab)^{t-1}=(a-b)(a^{t-1}+b^{t-1})\leq 2(a-b)\leq 2(a-b)^{\varepsilon}$,
and for each $t\in[\varepsilon,1)$, by the concavity of the function
$x^t$:
$a^t-b^t\leq (a-b)^t\leq (a-b)^{\varepsilon}$.}
This gives,
by considering any $t$ and $r_1,\ldots,r_k$ in the suprema for $x$, that
$\phi_k(y)\geq\phi_k(x)-2kd_R(x,y)^{\varepsilon}$
(using that the denominator in \rf{5ok12a} is at least 1).
So $\phi_k$ is continuous in the $d_R$-topology on $B(H)$.

Define for each $z\in B(X)$:
\dez{
U_z:=\{x\in B(H)\mid \phi_{f(z)}(x-z)<\varepsilon\}.
}
So $U_z$ is open in de $d_R$-topology.
Moreover, the $U_z$ for $z\in B(X)$ cover $B(H)$.
Indeed, for any $x\in B(H)$ there exists $z\in B(X)$ with
$\|x-z\|<\varepsilon^{1/\varepsilon}$.
Then $x\in U_z$, since $\phi_k(x-z)<\varepsilon$ for any $k$, which
follows from the following inequality.
Let $t\in[\varepsilon,2]$ and
$r_1,\ldots,r_k\in R$ be orthogonal and nonzero, for some $k\geq 1$.
Define $s_i:=r_i/\|r_i\|$ for each $i$.
So $s_1,\ldots,s_k$ are orthonormal.
Denote
$\rho:=\|(\|r_1\|^t,\ldots,\|r_k\|^t)\|_p$, with $p:=2/(2-t)$.
Then one has for any $y\in B(H)$, using the H\"older inequality,
and setting $q:=2/t$ (so that $p^{-1}+q^{-1}=1$):
\dyyz{
\sum_{i=1}^k |\langle r_i,y\rangle|^t
=
\sum_{i=1}^k\|r_i\|^t\cdot |\langle s_i,y\rangle|^t
\leq
(\sum_{i=1}^k\|r_i\|^{tp})^{1/p}
\cdot
(\sum_{i=1}^k|\langle s_i,y\rangle|^{tq})^{1/q}
=
\rho(\sum_{i=1}^k\langle s_i,y\rangle^2)^{1/q}
\leq
\rho\|y\|^{2/q}
=
\rho\|y\|^t
\leq
(1+\rho)\|y\|^{\varepsilon}.
}
So $\phi_{f(z)}(x-z)\leq\|x-z\|^{\varepsilon}<\varepsilon$, and hence $x\in U_z$.

As $(B(H),\delta_R)$ is compact by (iv), there is a finite set $Z\subseteq X$
such that voor each $x\in X$ there exist $z\in Z$ and $\pi\in G$ such that
$x\in U_{z^{\pi}}$.
This gives (iii).
\bx

\sectz{Applications}

Since $R$ spans $X$, $X$ is fully determined by the positive semidefinite $R\times R$ matrix
giving the inner products of pairs from $R$.
Then $G$ is given by a group of permutations of $R$ that
leave the matrix invariant.
It is convenient to realize that $R$ is weakly regular if (but not only if)
the orbit space $R^k/G$ is compact for each $k$.

\medskip
\noindent
{\bf 1. Szemer\'edi's regularity lemma}
[7].
Let $R$ be the collection of sets $I\times J$, with $I$ and $J$ each being a union of
finitely many subintervals of $[0,1]$,
with inner product equal to the measure of the intersection.
Let $G$ be the group of permutations of the intervals of any partition of $[0,1]$ into
intervals.
Then $G$ acts on $R$.

Let $\Pi$ be the collection of partitions of $[0,1]$ into finitely many
sets, each being a union of finitely many intervals.
For $P,Q\in\Pi$, $P\leq Q$ if and only if $P$ is a refinement of $Q$.
This gives a lattice; let $\wedge$ be the meet.

For any $P\in\Pi$, let $L_P$ be subspace of $X$ spanned by the elements $I\times J$
with $I,J\in P$.
For any $x\in X$, let $x_P$ be the orthogonal projection of $x$ onto $L_P$.

\lemma{4ok12a}{
For each $x\in X$ and $\varepsilon>0$ there exists $t_{\varepsilon,x}$
such that for each $N\in\Pi$ there is a $P\geq N$ such that
$\|x_N-x_P\|<\varepsilon$ and $|P|\leq t_{\varepsilon,x}$.
}

\pf
Let $Y$ be the set of those $x$ for which the statement holds for all $\varepsilon>0$.
Then $Y$ is a linear space.
Indeed, if $x\in Y$ and $\lambda\neq 0$ then $\lambda x\in Y$, as
we can take $t_{\varepsilon,\lambda x}:=t_{|\lambda^{-1}|\varepsilon,x}$.
If $x,y\in Y$ then $x+y\in Y$, as we can take
$t_{\varepsilon,x+y}:=t_{\varepsilon/2,x}t_{\varepsilon/2,y}$,
since if $\|x_N-x_P\|<\varepsilon/2$ and
$\|y_N-y_Q|<\varepsilon/2$ for some $P,Q\geq N$, then
$\|(x+y)_N-x_P-y_Q\|<\varepsilon$, hence $\|(x+y)_N-(x+y)_{P\wedge Q}\|\leq\varepsilon$,
since $x_P+y_Q\in L_{P\wedge Q}$ and $((x+y)_N)_{P\wedge Q}=(x+y)_{P\wedge Q}$ (since
$L_{P\wedge Q}\subseteq L_N$).
Note that $|P\wedge Q|\leq |P||Q|$.

So $Y$ is a linear space, and hence it suffices to show that $R\subseteq Y$.
Let $x\in R$ and $\varepsilon>0$.
We claim that $t_{\varepsilon,x}:=(1+2/\varepsilon)^2$ will do.
Indeed, let $N\in\Pi$.
Then
\dez{
x_N=\sum_{I,J\in N}\alpha_I\beta_J(I\times J)
}
for some $\alpha,\beta:N\to[0,1]$.
Let $\alpha'$ and $\beta'$ be obtained from $\alpha$ and $\beta$ by rounding down the values
to an integer multiple of $\varepsilon/2$.
Let $P\geq N$ be such that two classes $I$ and $J$ of $N$ are contained in the same class of $P$
if and only if $\alpha'_I=\alpha'_J$ and $\beta'_I=\beta'_J$.
As the pairs $(\alpha'_I,\beta'_I)$ take at most $(1+2/\varepsilon)^2$ different values,
we have $|P|\leq (1+2\varepsilon^{-1})^2$.
Define
\dez{
y:=\sum_{I,J\in N}\alpha'_I\beta'_J(I\times J).
}
Then $y\in L_P$.
Hence, since $x_P=(x_N)_P$ (as $L_P\subseteq L_N$), implying that $x_P$ is the point on $L_P$
closest to $x_N$:
\dyyz{
\|x_N-x_P\|^2
\leq
\|x_N-y\|^2
\leq
\sum_{I,J\in N}(\alpha_I\beta_J-\alpha'_I\beta'_J)^2\mu(I\times J)
\leq
\varepsilon^2
\sum_{I,J\in N}\mu(I\times J)
=
\varepsilon^2.
}
Here $\mu(I\times J)$ is the measure of $I\times J$.
\bx

Call a collection $P$ of sets {\em balanced} if all sets in $P$ have the same cardinality.
Call a partition $P$ of a finite set $V$ {\em $\varepsilon$-balanced} if
$P\setminus P'$ is balanced for some $P'\subseteq P$ with
$|\bigcup P'|\leq\varepsilon|V|$.

\lemma{12se12a}{
Let $\varepsilon>0$.
Then each partition $P$ of a finite set $V$ has an $\varepsilon$-balanced
refinement $Q$ with $|Q|\leq (1+1/\varepsilon)|P|$.
}

\pf
Define $t:=\varepsilon|V|/|P|$.
Split each class of $P$ into classes, each of size $\lceil t\rceil$,
except for at most one of size less than $t$.
This gives $Q$.
Then $|Q|\leq |P|+|V|/t=(1+1/\varepsilon)|P|$.
Moreover, the union of the classes of $Q$ of size less than $t$ has
size at most $|P|t=\varepsilon|V|$.
So $Q$ is $\varepsilon$-balanced.
\bx

Given a graph $H=(V,E)$ and $C,D\subseteq V$, then $e(C,D)$ is the number of
adjacent pairs of vertices in $C\times D$.
If $C,D\neq\emptyset$, let $d(C,D):=e(C,D)/|C||D|$.

\thmnmz{Szemer\'edi's regularity lemma}{
For each $\varepsilon>0$ and $p\in\oN$ there exists $k_{p,\varepsilon}\in\oN$
such that for each graph $H=(V,E)$ and each partition $P$ of $V$ with $|P|=p$
there is an $\varepsilon$-balanced refinement $Q$ of $P$ with $|Q|\leq k_{p,\varepsilon}$
and
\de{12se12c}{
\sum_{A,B\in Q}
\max_{\emptyset\neq C\subseteq A\atop \emptyset\neq D\subseteq B}
(|C||D|\cdot|d(C,D)-d(A\times B)|)^2
<
\varepsilon|V|^2.
}
}

\pf
Let $R$ and $G$ be as above.
It is easy to check that $R^k/G$ is compact for each $k$,
hence $(R,G)$ is very weakly regular.
So, by Theorem \ref{3ok12b}, $(R,G)$ is very strongly regular.

Fix $\varepsilon>0$ and $p\in\oN$.
For each $x\in X$, define
$f(x):=((1+1/\varepsilon)pt_{\varepsilon/4,x})^2$,
where $t_{\varepsilon/4,x}$ is as given in Lemma \ref{4ok12a}.

By the very strong regularity of $(R,G)$, there exists a finite set $Z\subseteq X$ such that
for each $x\in B(X)$ there exist $z\in Z$ and $\pi\in G$ satisfying
\vspace{-4mm}
\dy{25se12e}{
$\dps\sum_{j=1}^{f(z)}\langle r_j,x-z^{\pi}\rangle^2<\varepsilon^2/16$
for all orthogonal $r_1,\ldots,r_{f(z)}\in R$.
}
\vspace{-4mm}
Let $k_{p,\varepsilon}:=\max\{f(z)\mid z\in Z\}$.
We show that $k_{p,\varepsilon}$ is as required.

Let $H=([n],E)$ be a graph.
Let $N$ be the partition of $[0,1]$ into $n$ equal consecutive intervals $I_1,\ldots,I_n$,
and let $x:=\sum_{i,j\in [n]\text{ adjacent}}I_i\times I_j$
(the corresponding {\em graphon}).

By the above there exists a $z\in Z$ and a $\pi\in G$ satisfying
\rf{25se12e}.
By Lemma \ref{4ok12a}, there is a partition $U\in\Pi$ with $U\geq N$
such that $|U|\leq t_{\varepsilon/4,z}$ and $\|z_N-z_U\|\leq\varepsilon/4$.
Let $S:=P\wedge U$.
So $|S|\leq |P||U|\leq pt_{\varepsilon/4,z}$.
By Lemma \ref{12se12a}, there is an $\varepsilon$-balanced refinement $Q$ of $S$
with $N\leq Q\leq S$ and $|Q|\leq (1+1/\varepsilon)|S|\leq\sqrt{f(z)}\leq k_{p,\varepsilon}$.
We show that this $Q$ gives the partition of the theorem.

For each $A,B\in Q$, choose $r\in R$ with $r\subseteq A\times B$,
such that $r\in L_N$ and
such that $|\langle r,x-z_Q\rangle|$ is maximized.
This implies for each $r'\in R$ with $r'\subseteq A\times B$ and $r'\in L_N$:
\dyy{4ok12c}{
|\langle r',x-x_Q\rangle|
\leq
|\langle r',x-z_Q\rangle|
+
|\langle r',x_Q-z_Q\rangle|
\leq
|\langle r',x-z_Q\rangle|
+
|\langle A\times B,x_Q-z_Q\rangle|
=
|\langle r',x-z_Q\rangle|
+
|\langle A\times B,x-z_Q\rangle|
\leq
2|\langle r,x-z_Q\rangle|.
}

Let $r_1,\ldots,r_t$ be the chosen elements.
So $t=|Q|^2\leq f(z)$.
Hence, noting that $\langle r_i,z\rangle=\langle r_i,z_N\rangle$, since $r_i\in L_N$,
\dyy{4ok12d}{
(\sum_{i=1}^t\langle r_i,x-z_Q\rangle^2)^{1/2}
\leq
(\sum_{i=1}^t\langle r_i,x-z_N\rangle^2)^{1/2}
+
\|z_N-z_Q\|
\leq
(\sum_{i=1}^t\langle r_i,x-z\rangle^2)^{1/2}
+
\varepsilon/4
\leq
\varepsilon/2.
}
For the graph $H$, \rf{4ok12c} and \rf{4ok12d} give \rf{12se12c}.
\bx

To interpret \rf{12se12c}, for $A,B\in Q$, let $m_{A,B}$ denote the maximum
described in \rf{12se12c}.
Let $Q'$ be such that $Q\setminus Q'$ is balanced and $|\bigcup Q'|\leq\varepsilon|V|$.
Set $Q'':=Q\setminus Q'$, and let $Z$ be the collection of pairs $(A,B)\in Q''\times Q''$
with $m_{A,B}\geq\sqrt{\varepsilon}|A||B|$.
Then \rf{12se12c} implies
\dez{
\sum_{(A,B)\in Z}|A||B|
\leq
\sum_{(A,B)\in Z}\varepsilon^{-1/2}m_{A,B}
\leq
\sqrt{\varepsilon}|V|^2.
}
Moreover, as $|\bigcup Q'|<\varepsilon|V|$,
\dyyz{
\sum_{A,B\in Q''}|A||B|\geq \sum_{A,B\in Q}|A||B|-2\varepsilon|V|^2=(1-2\varepsilon)|V|^2.
}
Hence, assuming $\varepsilon<1/4$, $|Z|\leq \sqrt{\varepsilon}(1-2\varepsilon)^{-1}|Q''|^2<2\sqrt{\varepsilon}|Q''|^2$.
For each $(A,B)\in(Q''\times Q'')\setminus Z$ one has $m_{A,B}<\sqrt{\varepsilon}|A||B|$,
implying that for each rectangle $R\subseteq A\times B$ with
$|R|/|A\times B|\geq\sqrt[4]{\varepsilon}$ one has $|d(R)-d(A\times B)|<\sqrt[4]{\varepsilon}$.
In other words, $A\times B$ is $\sqrt[4]{\varepsilon}$-regular.

\medskip
\noindent
{\bf 2. ``Interval regularity''.}
Let $R$ be the collection of sets $I\times J$, with $I$ and $J$ subintervals of $[0,1]$,
with inner product given by the measure of the intersection.
Then Theorem \ref{3ok12b} gives an ``interval regularity theorem'' for graphs
(it can also be proved with Szemer\'edi's classical combinatorial method):

\thmz{
For each $\varepsilon>0$ and $p\in\oN$ there exists $k_{p,\varepsilon}\in\oN$ such that
for each $n$, each graph $H=([n],E)$ and each partition $P$ of $[n]$ into intervals
with $|P|\leq p$, $P$ has a refinement to a partition $Q$ into at most $k_{p,\varepsilon}$
intervals such that all intervals in $Q$ have the same size except for some of them
covering $\leq\varepsilon n$ vertices and such that
\dez{
\sum_{A,B\in Q}\max_{I\subseteq A, J\subseteq B\atop I,J\text{\rm ~intervals}}|I||J||d(I,J)-d(A,B)|
<
\varepsilon n^2.
}
{\rm Here $d(I,J)$ and $d(A,B)$ are the densities of the corresponding subgraphs of $H$.}
}

This can be derived similarly as (in fact, easier than) Szemer\'edi's regularity lemma above.

\medskip
\noindent
{\bf 3. Polynomial approximation.}
Let $k\leq n$.
Each polynomial $p\in\oR[x_1,\ldots,x_n]$ can be uniquely written as
$p=\sum_{\mu}\mu p_{\mu}$, where $\mu$ ranges over the set $M$ of all monomials in $\oR[x_1,\ldots,x_k]$
and where $p_{\mu}\in\oR[x_{k+1},\ldots,x_n]$.
If $p$ is homogeneous of degree $d$, we say
that $p$ is {\em $\varepsilon$-concentrated on the first $k$ variables} if
\dez{
\sum_{\mu\in M\atop\deg(\mu)<d}\max_{x\in\oR^{n-k}\atop\|x\|=1}p_{\mu}(x)^2
\leq
\varepsilon\|p\|^2,
}
where $\|p\|$ is the square root of the sum of the squares of the coefficients of $p$.

\thmz{
For each $\varepsilon>0$ and $d\in\oN$ there exists $k_{d,\varepsilon}$ such that for each $n$,
each homogeneous polynomial of degree $d$ in $n$ variables is $\varepsilon$-concentrated on the first $k$
variables
\textit{\textbf{after some orthogonal transformation of $\oR^n$}},
for some $k\leq k_{d,\varepsilon}$.
}

\medskip
This can be derived by setting $R$ to be the set of all polynomials
$(a\T x)^d$, with $a\in\oR^n$ and $\|a\|=1$ for some $n$ (setting
$x=(x_1,x_2,\ldots)$), taking the inner product of
$(a\T x)^d$ and $(b\T x)^d$ equal to $(a\T b)^d$.
(This corollary strengthens a `weak regularity' result of
Fernandez de la Vega, Kannan, Karpinski, and Vempala [2].)
For details, we refer to [6].

\section*{References}\label{REF}
{\small
\begin{itemize}{}{
\setlength{\labelwidth}{4mm}
\setlength{\parsep}{0mm}
\setlength{\itemsep}{1mm}
\setlength{\leftmargin}{5mm}
\setlength{\labelsep}{1mm}
}
\item[\mbox{\rm[1]}] J. Dugundji, 
{\em Topology},
Allyn and Bacon, Boston, 1966.

\item[\mbox{\rm[2]}] W. Fernandez de la Vega, R. Kannan, M. Karpinski, S. Vempala, 
Tensor decomposition and approximation schemes for constraint satisfaction problems,
in: {\em Proceedings of the 37th Annual {ACM} Symposium on Theory of Computing}
({STOC}'05),
pp. 747--754,
{ACM}, New York, 2005.

\item[\mbox{\rm[3]}] L. Lov\'asz, B. Szegedy, 
Limits of dense graph sequences,
{\em Journal of Combinatorial Theory, Series B} 96 (2006) 933--957.

\item[\mbox{\rm[4]}] L. Lov\'asz, B. Szegedy, 
Szemer\'edi's lemma for the analyst,
{\em Geometric and Functional Analysis} 17 (2007) 252--270.

\item[\mbox{\rm[5]}] G. Regts, A. Schrijver, 
Compact orbit spaces in Hilbert spaces and limits of edge-colouring models,
preprint, 2012.
ArXiv \url{http://arxiv.org/abs/1210.2204}

\item[\mbox{\rm[6]}] A. Schrijver, 
Low rank approximation of polynomials,
preprint, 2012.\\
\url{http://www.cwi.nl/~lex/lrap.pdf}

\item[\mbox{\rm[7]}] E. Szemer\'edi, 
Regular partitions of graphs,
in: {\em Probl\`emes combinatoires et th\'eorie des graphes}
(Proceedings Colloque International C.N.R.S., Paris-Orsay, 1976;
J.-C. Bermond, J.-C. Four\-nier, M. Las Vergnas, D. Sotteau, eds.)
[Colloques Internationaux du Centre National de la Recherche
Scientifique N$^o$ 260],
\'Editions du Centre National de la Recherche Scientifique,
Paris, 1978, pp. 399--401.

\end{itemize}
}

\end{document}